\documentclass[11pt,a4paper]{article}

\usepackage{inputenc}
\usepackage{amsmath}
\usepackage{bm}
\usepackage{bbold}
\usepackage{amsthm}
\usepackage{enumerate}

\usepackage{algorithm}

\usepackage[top=1.0in, left=1.0in, right=1.0in, bottom=1.0in]{geometry}

\usepackage{tikz}\tikzset{x=1cm,y=1cm,z=1cm}

\usepackage{pgfplots}\pgfplotsset{compat=1.16}

\usepackage[hyphens]{url}

\usepackage{hyperref}
\usepackage{breakurl}

\title{Application of log-Chebyshev approximation and tropical algebra to multicriteria problems of pairwise comparisons}

\author{N. Krivulin\thanks{Faculty of Mathematics and Mechanics, St.~Petersburg State University, 28 Universitetsky Ave., St.~Petersburg, 198504, Russia; 
nkk@math.spbu.ru.}
}

\date{}

\newtheorem{theorem}{Theorem}
\newtheorem{lemma}[theorem]{Lemma}

\theoremstyle{definition}

\begin{document}

\maketitle

\begin{abstract}
We consider multicriteria problems of evaluating absolute ratings (scores, priorities, weights) of given alternatives for making decisions, which are compared in pairs under several criteria. Given matrices of pairwise comparisons of alternatives for each criterion and a matrix of pairwise comparisons of the criteria, the aim is to calculate a vector of individual ratings of alternatives. We formulate the problem as the Chebyshev approximation of matrices on the logarithmic scale by a common consistent matrix (a symmetrically reciprocal matrix of unit rank). We rearrange the approximation problem as a multi-objective optimization problem of finding a vector that determines the consistent matrix and hence yields a vector of ratings in question. The problem is then transformed into a series of optimization problems in the framework of tropical algebra, which focuses on the theory and application of algebraic systems with idempotent operations. To solve the optimization problems, we apply methods and results of tropical optimization, which yield analytical solutions in a form ready for further analysis and straightforward computation. To illustrate the technique implemented, we give a numerical example of solving a known problem, and compare the obtained solution with results provided by classical methods of analytic hierarchy process and weighted geometrical means.
\\

\\

\end{abstract}

\section{Introduction}
\label{S-I}

Multicriteria problems of evaluating alternatives based on pairwise comparisons are encountered in many areas of practice and constitute an important class of problems that continue to attract considerable attention from decision-making specialists \cite{Saaty1990Analytic,Gavalec2015Decision,Ramik2020Pairwise,Greco2025Fifty}. These problems involve finding absolute ratings (scores, priorities, weights) of alternatives for subsequent selection of the best alternative for making a decision, based on results of pairwise comparisons of the alternatives according to given criteria, as well as results of pairwise comparisons of the criteria themselves. The results can serve to order (rank) the alternatives according to the magnitude of their ratings.

To solve pairwise comparison problems, various tools are used \cite{Choo2004Common}, including heuristic procedures that do not guarantee the optimality of the obtained result but usually provide acceptable accuracy for the application, as well as formal methods that lead to mathematically sound results but may require solving difficult computational problems. The most well-known examples are the heuristic method of analytic hierarchy process \cite{Saaty1977Scaling,Saaty1984Comparison,Saaty1990Analytic,Saaty2013Measurement} and the formal method of weighted geometric means \cite{Narasimhan1982Geometric,Belton1983Shortcoming,Crawford1985Note}, which allow obtaining results with low computational complexity.

At the same time, it is known that existing methods for solving the problems can yield different and even contradictory results in practice \cite{Saaty1984Comparison,Barzilai1987Consistent,Ishizaka2006How,Tran2013Pairwise,Mazurek2021Numerical}. The possibility of contradictory results complicates the use of these methods as effective decision support tools, and thus the need for further expansion of the available tools and comparative analysis of the resulting solutions remains relevant. When the results obtained using different methods contradict each other, their application for decision making seems insufficiently justified. On the other hand, agreement or similarity in the obtained results may serve as an argument in favor of any of them as a solution close to optimal.

One of the methods that offer alternative approaches to solving multicriteria pairwise comparison problems is a method based on the log-Chebyshev approximation of pairwise comparison matrices, proposed and studied in a series of work \cite{Krivulin2015Computation,Krivulin2015Rating,Krivulin2016Using,Krivulin2017Application,Krivulin2018Methods,Krivulin2019Methods,Krivulin2019Tropical,Krivulin2020Using,Krivulin2021Algebraic,Krivulin2024Application}. The solution offered by the method reduces the problem to solving a number of optimization problems in terms of tropical algebra which deals with theory and application of algebraic systems with idempotent operations \cite{Kolokoltsov1997Idempotent,Golan2003Semirings,Heidergott2006Maxplus,Maclagan2015Introduction,Kenoufi2025Idempotent}. This approach leads to an analytical solution to the problem, which represents all solutions in a compact vector form, suitable for analysis and calculations.

In this paper we consider multicriteria problems of evaluating absolute ratings of given alternatives for making decisions, which are compared in pairs under several criteria. Given matrices of pairwise comparisons of alternatives for each criterion and a matrix of pairwise comparisons of the criteria, the aim is to calculate a vector of individual ratings of alternatives. We formulate the problem as the Chebyshev approximation of matrices on the logarithmic scale by a common consistent matrix (a symmetrically reciprocal matrix of unit rank). We rearrange the approximation problem as a multi-objective optimization problem of finding a vector that determines the consistent matrix and hence yields a vector of ratings in question. The problem is then transformed into a series of optimization problems in the framework of tropical algebra, which focuses on the theory and applications of algebraic systems with idempotent operations. To solve the optimization problems, we apply methods and results of tropical optimization, which yield analytical solutions in a form ready for further analysis and straightforward computation. To illustrate the technique implemented, we give a numerical example of solving a known problem, and compare the obtained solution with results provided by classical methods of analytic hierarchy process and weighted geometrical means.

The rest of the paper is organized as follows. In Section~\ref{S-SCPCP}, we consider a pairwise comparison problem with one criterion and discuss its solutions, including the solution on the basis of log-Chebyshev approximation. Section~\ref{S-MPCP} describes a general multicriteria problem and outline three methods of solution: analytic hierarchy process, weighted geometric means and log-Chebyshev approximation. In Section~\ref{S-LCAUTA}, we show how the multicriteria pairwise comparison problem can be solved using log-Chebyshev approximation in terms of tropical algebra. A numerical example is given in Section~\ref{S-SSE}.

\section{Single Criterion Pairwise Comparison Problem}
\label{S-SCPCP}

Let $\mathcal{A}_{1},\ldots,\mathcal{A}_{n}$ be $n$ alternatives available for making a decision. Suppose that the alternatives are compared pairwise, resulting in an $(n\times n)$-matrix $\bm{A}=(a_{ij})$ of pairwise comparisons, where the element $a_{ij}>0$ indicates how many times alternative $\mathcal{A}_{i}$ is considered more preferable than alternative $\mathcal{A}_{j}$. It is assumed that the entries in the pairwise comparison matrix $\bm{A}$ satisfies the condition $a_{ij}=1/a_{ji}$ for all $i,j=1,\ldots,n$ (and hence $a_{ii}=1$), which makes the matrix $\bm{A}$ symmetrically reciprocal. This condition states that if alternative $\mathcal{A}_{i}$ is $a_{ij}$ times better than $\mathcal{A}_{j}$, then alternative $\mathcal{A}_{j}$ must be $1/a_{ij}$ times better ($a_{ij}$ times worse) than $\mathcal{A}_{i}$. Given a pairwise comparison matrix $\bm{A}$ representing the results of the relative evaluation of one alternative with respect to another, the problem of interest is to calculate the absolute ratings of the alternatives in the form of an $n$-vector $\bm{x}=(x_{j})$, where $x_{j}>0$ indicates the rating of alternative $\mathcal{A}_{j}$.

A pairwise comparison matrix $\bm{A}$ is called consistent if the transitivity property $a_{ij}=a_{ik}a_{kj}$ holds for all $i,j,k=1,\ldots,n$. This condition corresponds to the natural transitivity of judgments, which implies that if alternative $\mathcal{A}_{i}$ is considered $a_{ik}$ times better than $\mathcal{A}_{k}$, and alternative $\mathcal{A}_{k}$ is $a_{kj}$ times better than $\mathcal{A}_{j}$, then alternative $\mathcal{A}_{i}$ must be $a_{ik}a_{kj}$ times better than $\mathcal{A}_{j}$. The transitivity property of a consistent matrix $\bm{A}$ yields the existence of a positive vector $\bm{x}=(x_{j})$ whose entries determine the entries of $\bm{A}$ by the relation $a_{ij}=x_{i}/x_{j}$ for all $i,j$, which means that the matrix $\bm{A}$ is of rank one.

Pairwise comparison matrices encountered in real-world problems are usually not consistent, which makes the problem of deriving absolute ratings of alternatives from the results of pairwise comparisons (the pairwise comparison problem) challenging. Available methods for solving this problem include various heuristic methods, which in many cases in practice give acceptable results, as well as approximation methods, which provide mathematically sound optimal solutions.

\subsection{Solution Approaches to Pairwise Comparison Problem}
\label{S-SAPCP}

In this section, we consider solution methods, which focus on the derivation of a consistent matrix $\bm{X}=(x_{i}/x_{j})$ that is close to a given inconsistent pairwise comparison matrix $\bm{A}$ in some sense and thus determines an approximate solution to the pairwise comparison problem. Since in the context of this problem, the main goal is to find the vector of ratings $\bm{x}=(x_{j})$ rather than its associated matrix $\bm{X}$, the methods are commonly represented as a procedure of direct evaluation of this vector. 

Heuristic methods are generally based on aggregating the columns of a pairwise comparison matrix $\bm{A}$ to obtain a solution by calculating a weighted sum of these columns \cite{Choo2004Common}. In the most widely used principal eigenvector method \cite{Saaty1977Scaling,Saaty1990Analytic,Saaty2013Measurement}, the solution vector is defined as the sum of the columns taken with weights proportional to the entries of $\bm{x}$. In this case, the vector $\bm{x}$ satisfies the equation $\bm{A}\bm{x}=\lambda\bm{x}$ where $1/\lambda$ is the proportionality coefficient, and it is obtained as the principal (Perron) eigenvector of the matrix $\bm{A}$, which corresponds to the eigenvalue with the largest absolute value.

Approximation methods reduce the problem to finding a consistent matrix $\bm{X}$ that approximates a given pairwise comparison matrix $\bm{A}$ in the sense of some distance function as the approximation error \cite{Choo2004Common}. The approximation problem is solved by minimizing the distance between the matrices $\bm{A}$ and $\bm{X}$, which provides a formal justification for the solution obtained. The positive vector $\bm{x}$, which defines the approximating consistent matrix, is taken as the vector of absolute ratings of the alternatives.

If the error is measured on a standard linear scale, the approximation approach typically leads to complex multiextremal optimization problems that are difficult to solve \cite{Saaty1984Comparison}. In contrast, using a logarithmic scale (with a logarithm to a base greater than 1) makes the approximation problem more feasible and even allows for the analytical derivation of a solution in an exact, explicit form.

A widely used approximation technique for solving the pairwise comparison problem is the geometric mean method, which measures the error between the matrices $\bm{A}=(a_{ij})$ and $\bm{X}=(x_{i}/x_{j})$ using the Euclidean metric on a logarithmic scale (log-Euclidean approximation) \cite{Narasimhan1982Geometric,Crawford1985Note,Barzilai1987Consistent}. The method consists of finding a positive vector $\bm{x}=(x_{j})$, which solves the minimization problem
\begin{equation*}
\begin{aligned}
\min_{\bm{x}>\bm{0}}
&&&
\sum_{1\leq i,j\leq n}\left(\log a_{ij}-\log\frac{x_{i}}{x_{j}}\right)^{2}.
\end{aligned}
\end{equation*}

The usual solution technique, which applies the first derivative test to find the roots of the partial derivatives of the objective function with respect to all $x_{i}$, yields a solution vector $\bm{x}$ with entries specified in the parametric form
\begin{equation*}
x_{i}
=
\left(
\prod_{j=1}^{n}a_{ij}
\right)^{1/n}
u,
\qquad
u>0,
\qquad
i=1,\ldots,n.
\end{equation*}

This vector is normalized (for example, by dividing by the sum of the entries), and the result is taken as the optimal solution to the pairwise comparison problem.

Another approximation technique developed in \cite{Krivulin2015Computation,Krivulin2015Rating,Krivulin2016Using} and other works proposes a solution using the Chebyshev distance on a logarithmic scale (the log-Chebyshev approximation). The method aims at finding positive vectors $\bm{x}$ that solve the problem
\begin{equation}
\begin{aligned}
\min_{\bm{x}>\bm{0}}
&&&
\max_{1\leq i,j\leq n}\left|\log a_{ij}-\log\frac{x_{i}}{x_{j}}\right|.
\end{aligned}
\label{P-minx_maxijlogaijlogxixj}
\end{equation}

The solution of the problem using optimization methods in terms of conventional mathematics is difficult due to the complicated nonsmooth multiextremal objective function involved. Below we show how this problem can be solved in the framework of tropical optimization in direct analytical form.

We note that the principal eigenvector and geometric mean methods lead to unique solutions and offer efficient computational procedures for calculating the result. At the same time, the solution obtained using the log-Chebyshev approximation may be nonunique and therefore requires further analysis. However, given the approximate character of the pairwise comparison model, which usually leads to inconsistent results of comparisons, having multiple solutions to the problem seems quite reasonable and even allows one to choose a solution that satisfies some additional constraints.

\subsection{Log-Chebyshev Approximation of Pairwise Comparison Matrices}

Consider the log-Chebyshev approximation problem given by \eqref{P-minx_maxijlogaijlogxixj}. Since the logarithm monotonically increases, and the matrices $\bm{A}$ and $\bm{X}$ are reciprocal and have only positive entries, the objective function in \eqref{P-minx_maxijlogaijlogxixj} can be represented as (see \cite{Krivulin2017Application,Krivulin2018Methods,Krivulin2019Methods})
\begin{equation*}
\max_{1\leq i,j\leq n}\left|\log a_{ij}-\log\frac{x_{i}}{x_{j}}\right|
=
\log\max_{1\leq i,j\leq n}\frac{a_{ij}x_{j}}{x_{i}}.
\end{equation*}

The minimization of the logarithmic function on the right-hand side is equivalent to minimization of its argument, and thus we drop the logarithm from the objective function to reduce problem \eqref{P-minx_maxijlogaijlogxixj} to the problem
\begin{equation}
\begin{aligned}
\min_{\bm{x}>\bm{0}}
&&&
\max_{1\leq i,j\leq n}\frac{a_{ij}x_{j}}{x_{i}}.
\end{aligned}
\label{P-minx_maxijaijxixj}
\end{equation}

It is not difficult to verify (see, e.g., \cite{Krivulin2015Computation,Krivulin2020Using}) that in the case of a reciprocal matrix $\bm{A}$ and consistent rank one matrix $\bm{X}$, the objective function in \eqref{P-minx_maxijaijxixj} satisfies the equality
\begin{equation*}  
\max_{1\leq i,j\leq n}\frac{a_{ij}x_{j}}{x_{i}}
=
\max_{1\leq i,j\leq n}\frac{|a_{ij}-x_{i}/x_{j}|}{a_{ij}}
+
1,
\end{equation*}
where the first term on the right-hand side takes the form of the maximum relative error. 

It follows from the equality that the objective function reaches its minimum simultaneously with this error. As a result, the solution of problem \eqref{P-minx_maxijaijxixj} is equivalent to the approximation in the sense of the maximum relative error as well. 

As it is shown later, the problem of log-Chebyshev approximation at \eqref{P-minx_maxijaijxixj} can be directly solved in terms of tropical algebra to give the result in a compact vector form.

\subsection{Best and Worst Differentiating Solutions}

Suppose that the solution of problem \eqref{P-minx_maxijaijxixj} yields a vector $\bm{x}=(x_{j})$ that is uniquely determined (up to a positive factor). In this case, the vector $\bm{x}$ is taken as the vector of absolute ratings of alternatives under discussion. 

If the solution of \eqref{P-minx_maxijaijxixj} is not unique, consider the set of all solution vectors
\begin{equation}
X
=
\arg\min_{\bm{x}>\bm{0}}
\max_{1\leq i,j\leq n}\frac{a_{ij}x_{j}}{x_{i}}.
\label{E-X}
\end{equation}

Since nonunique solutions make it difficult to derive an optimal decision in the pairwise comparison problem, further investigation is necessary to characterize the result by one or two vectors that are reasonably representative of the set of solutions. To overcome this difficulty, an analytical approach developed in \cite{Krivulin2017Application,Krivulin2018Methods,Krivulin2019Tropical} suggests to reduce the entire set of solutions to two vectors that can be considered as the ``best'' and ``worst'' solutions in the sense of better differentiation of the alternatives with the highest and lowest ratings. As the best solution, the approach chooses a vector of ratings with the maximal ratio between its maximum and minimum entries, and as the worst a vector with the minimal ratio between these entries.

The ratio between the highest and lowest ratings given by a vector $\bm{x}$ takes the form of the multiplicative Hilbert (span, range) seminorm of $\bm{x}$, which is written as
\begin{equation}
\max_{1\leq i\leq n}x_{i}
\Big/
\min_{1\leq j\leq n}x_{j}
=
\max_{1\leq i\leq n}x_{i}
\times
\max_{1\leq j\leq n}x_{j}^{-1}.
\label{E-maxximinxi}
\end{equation}

The best and worst differentiating solutions are found by maximizing and minimizing the Hilbert seminorm over all vectors $\bm{x}\in X$, which leads to the optimization problems
\begin{equation*}
\begin{aligned}
\max_{\bm{x}\in X}
&&&
\max_{1\leq i\leq n}x_{i}\times\max_{1\leq j\leq n}x_{j}^{-1};
\end{aligned}
\qquad\qquad
\begin{aligned}
\min_{\bm{x}\in X}
&&&
\max_{1\leq i\leq n}x_{i}\times\max_{1\leq j\leq n}x_{j}^{-1}.
\end{aligned}
\end{equation*}

One of the shortcomings of this approach is that the best and worst differentiating solutions obtained in this way may not be unique as well. To handle this possible issue, an improved procedure of finding the best and worst differentiating solutions is proposed in \cite{Krivulin2024Application}, which results in a reduced set of the best solutions and a unique worst solution. 

The procedure takes into account that the solutions of problem \eqref{P-minx_maxijaijxixj} are invariant under multiplication by a positive factor, and hence concentrates only on normalized solutions obtained by dividing the solution vectors by their maximal entry. Adding the normalization condition leads to the optimization problems
\begin{equation}
\begin{aligned}
\max_{\bm{x}\in X}
&&&
\max_{1\leq i\leq n}x_{i}^{-1};
\\
\text{s.t.}
&&&
\max_{1\leq i\leq n}x_{i}=1;
\end{aligned}
\qquad\qquad
\begin{aligned}
\min_{\bm{x}\in X}
&&&
\max_{1\leq i\leq n}x_{i}^{-1};
\\
\text{s.t.}
&&&
\max_{1\leq i\leq n}x_{i}=1.
\end{aligned}
\label{P-maxxmaxi1xi-maxixi1-minxmaxi1xi-maxixi1}
\end{equation}

All normalized solution vectors that maximize or minimize the ratio between the highest and lowest ratings have two entries whose ratio is fixed: the maximum entry equal to 1, and the minimum entry less or equal to 1. In this case, the lower (higher) the other entries, the better (worse) the alternative with the maximum rating is distinguishable from the others. As a result, it makes sense to consider only the minimal normalized solution to the maximization problem and the maximal normalized solution to the minimization problem. As usual, a solution $\bm{x}_{0}$ is referred to as minimal (maximal) if for any solution $\bm{x}$ the componentwise inequality $\bm{x}_{0}\leq\bm{x}$ ($\bm{x}_{0}\geq\bm{x}$) holds.

\section{Multicriteria Pairwise Comparison Problem}
\label{S-MPCP}

Consider the problem of evaluating the ratings of alternatives, in which $n$ alternatives $\mathcal{A}_{1},\ldots,\mathcal{A}_{n}$ are compared pairwise according to $m$ criteria $\mathcal{C}_{1},\ldots,\mathcal{C}_{m}$ (the multicriteria pairwise comparison problem). Let $\bm{A}_{k}$ denote an $(n\times n)$-matrix of pairwise comparisons of the alternatives according to criterion $\mathcal{C}_{k}$. The criteria are also compared pairwise to form an $(m\times m)$-matrix $\bm{C}=(c_{kl})$ of pairwise comparisons for criteria, where $c_{kl}$ indicates how many times criterion $\mathcal{C}_{k}$ is more important for making decisions than $\mathcal{C}_{l}$. The problem is formulated to find the absolute individual rating of each alternative based on the pairwise comparison matrices $\bm{C}$ and $\bm{A}_{1},\ldots,\bm{A}_{m}$.

The most common approach to solving the problem is to use the heuristic method of Analytic Hierarchy Process (AHP) \cite{Saaty1990Analytic,Saaty2013Measurement}. The solution involves finding the normalized (with respect to the sum of entries) principal eigenvector $\bm{w}=(w_{k})$ for the matrix $\bm{C}$ and principal eigenvectors $\bm{x}_{1},\ldots,\bm{x}_{n}$ for the matrices $\bm{A}_{1},\ldots,\bm{A}_{m}$. Then the vector of absolute ratings of alternatives is calculated as the weighted sum
\begin{equation}
\bm{x}
=
\sum_{k=1}^{m}w_{k}\bm{x}_{k}.
\label{E-AHP}
\end{equation}

A direct approach to solving multicriteria pairwise comparison problems involves approximating the pairwise comparison matrices for all criteria with a common consistent matrix. The approximation problem reduces to a multiobjective optimization problem of minimizing a vector objective function of $m$ elements each of which determines the approximation error of the pairwise comparison matrix for a single criterion. The problem obtained is solved to find Pareto-optimal or, at least, weakly Pareto-optimal solutions using methods and techniques of multiobjective optimization \cite{Ehrgott2005Multicriteria,Luc2008Pareto,Benson2009Multiobjective,Nakayama2009Sequential}.

To solve multiobjective minimization problems, various scalarization methods are commonly used, which exploit the replacement of the vector objective function with a scalar one. The scalarization function is often chosen to be a weighted sum or weighted maximum of the elements of the vector objective function. Provided that all weights are nonzero, the solution obtained by minimizing a weighted sum (maximum) is known to be Pareto-optimal (weakly Pareto-optimal) (see, e.g., \cite{Nakayama2009Sequential}). If minimizing a weighted maximum gives a unique solution, then this solution is also Pareto optimal.

Approximation of pairwise comparison matrices with error measured on a linear scale leads to computationally complex problems even for a single criterion. At the same time, using a logarithmic scale to measure error allows for an explicit analytical solution to the pairwise comparison problem under consideration.

The solution to the multicriteria pairwise comparison problem using the log-Euclidean approximation is known as the method of Weighted Geometric Means (WGM). It applies scalarization of the objective function by the weighted sum
\begin{equation*}
\sum_{k=1}^{m}
w_{k}
\sum_{1\leq i,j\leq n}
\left(
\log a_{ij}^{(k)}-\log\frac{x_{i}}{x_{j}}
\right)^{2},
\end{equation*}
where the weights $w_{k}$ are obtained from the matrix $\bm{C}$ by the geometric mean method.

As a result of minimizing this sum (based on calculating the partial derivatives and equating them to zero), one obtains a vector $\bm{x}=(x_{i})$ with coordinates
\begin{equation}
x_{i}
=
\prod_{k=1}^{m}
\left(
\prod_{j=1}^{n}a_{ij}^{(k)}
\right)^{w_{k}/n}
u,
\qquad
i=1,\ldots,n;
\qquad
u>0.
\label{E-WGM}
\end{equation}

The resulting vector is a Pareto-optimal solution to the multiobjective optimization problem. This vector (usually normalized with respect to the sum of its elements) is taken as the WGM solution of the multicriteria pairwise comparison problem.

The solution of the multicriteria pairwise comparison problem using the $\log$-Chebyshev approximation (LCA) involves scalarization of the vector objective function into a weighted maximum of the approximation errors over all criteria in the form
\begin{equation*}
\max_{1\leq k\leq m}
w_{k}
\max_{1\leq i,j\leq n}
\frac{a_{ij}^{(k)}x_{j}}{x_{i}}
=
\max_{1\leq i,j\leq n}
\left(
\max_{1\leq k\leq m}
w_{k}a_{ij}^{(k)}
\right)
\frac{x_{j}}{x_{i}},
\end{equation*}
where the weights $w_{k}$ are obtained from the matrix $\bm{C}$ by solving the problem 
\begin{equation}
\begin{aligned}
\min_{\bm{w}>\bm{0}}
&&&
\max_{1\leq k,l\leq m}
\frac{c_{kl}w_{l}}{w_{k}}.
\end{aligned}
\label{P-LCA_w}
\end{equation}

Let $\bm{A}=(a_{ij})$ be an $(n\times n)$-matrix defined as the weighted maximum
\begin{equation}
\bm{A}
=
\max(w_{1}\bm{A}_{1},\ldots,w_{m}\bm{A}_{m}).
\label{E-Aeqw1A1wkAk}
\end{equation}

The multicriteria pairwise comparison problem under study is then formulated to find a positive vector $\bm{x}=(x_{j})$ that solves the minimax problem
\begin{equation}
\begin{aligned}
\min_{\bm{x}>\bm{0}}
&&&
\max_{1\leq i,j\leq n}
\frac{a_{ij}x_{i}}{x_{j}}.
\end{aligned}
\label{P-LCA_x}
\end{equation}

Provided that both problems \eqref{P-LCA_w} and \eqref{P-LCA_x} have unique solutions, then the LCA solution procedure for the multicriteria pairwise comparison problem consists of the following three main steps.
\begin{enumerate}
\item
Determination of the vector of weights $\bm{w}$ for the criteria by solving problem \eqref{P-LCA_w}.
\item
Calculation of the matrix $\bm{A}$ as the weighted maximum according to \eqref{E-Aeqw1A1wkAk}.
\item\label{Step3}
Determination of the vector of ratings $\bm{x}$ for the alternatives by solving problem \eqref{P-LCA_x}.
\end{enumerate}

If it turns out that the solution of problem \eqref{P-LCA_x} is not unique, then Step~\ref{Step3} is complemented by finding the best and worst vectors of ratings by solving the maximization and minimization problems at \eqref{P-maxxmaxi1xi-maxixi1-minxmaxi1xi-maxixi1} over the set $X$ given by \eqref{E-X}. 

If the solution of problem \eqref{P-LCA_w} is not unique, then the procedure is extended as suggested in \cite{Krivulin2019Methods,Krivulin2019Tropical} to include the following steps.
\begin{enumerate}
\item
Determination of the vectors of weights by solving problem \eqref{P-LCA_w}; solution of maximization and minimization problems at \eqref{P-maxxmaxi1xi-maxixi1-minxmaxi1xi-maxixi1} for these vectors to find the best differentiating vector of weights $\bm{w}=(w_{i})$ and the worst differentiating vector $\bm{v}=(v_{i})$.
\item
Calculation of the matrix of weighted maximum $\bm{P}=\max(w_{1}\bm{A}_{1},\ldots,w_{m}\bm{A}_{m})$.
\item
Determination of the vectors of ratings by solving problem \eqref{P-LCA_x} for $\bm{P}$; solution of maximization problem at \eqref{P-maxxmaxi1xi-maxixi1-minxmaxi1xi-maxixi1} to find the best differentiating vector of ratings $\bm{x}$.
\item
Calculation of the matrix of weighted maximum $\bm{R}=\max(v_{1}\bm{A}_{1},\ldots,v_{m}\bm{A}_{m})$.
\item
Determination of the vectors of ratings $\bm{y}$ by solving problem \eqref{P-LCA_x} for $\bm{R}$; solution of minimization problem at \eqref{P-maxxmaxi1xi-maxixi1-minxmaxi1xi-maxixi1} to find the worst differentiating vector of ratings $\bm{y}$.
\end{enumerate}

The extended procedure solves two subproblems, one with the best differentiating vector of weights $\bm{w}$, and the second with the worst vector $\bm{v}$. For the first subproblem, the best differentiating vector of ratings is taken as the best solution of the entire problem. The worst vector of ratings found in the second subproblem becomes the worst solution. 

We describe the procedure in more detail in the next section, where the solution of the optimization problems involved are given in terms of tropical algebra.

\section{Log-Chebyshev Approximation Using Tropical Algebra}
\label{S-LCAUTA}

We start this section with a short overview of the basic definitions and notation of tropical max-algebra, which are used to describe the solution of the pairwise comparison problem in the sequel. Further information on the theory, methods and applications of tropical mathematics can be found, for example, in \cite{Kolokoltsov1997Idempotent,Golan2003Semirings,Heidergott2006Maxplus,Maclagan2015Introduction,Kenoufi2025Idempotent}.

Furthermore, we present solutions of optimization problems that are formulated and solved in the framework of tropical algebra (tropical optimization problems). These results form the basis for the solution of multicriteria pairwise comparison problems, which is described as a computational procedure in the last part of the section.

\subsection{Elements of Tropical Algebra}

Tropical algebra is concerned with the theory and application of algebraic systems with idempotent operations. An operation is idempotent if its application to arguments of the same value yields this value as the result. For example, taking the maximum is idempotent since $\max(x,x)=x$, whereas arithmetic addition is not: $x+x=2x$.

Consider max-algebra, which is an algebraic system with an idempotent operation defined on the set of nonnegative real numbers $\mathbb{R}_{+}=\{x\in\mathbb{R}|x\geq0\}$. It is closed under addition denoted by $\oplus$ and performed for all $x,y\in\mathbb{R}_{+}$ as maximum: $x\oplus y=\max\{x,y\}$, and under multiplication defined as usual. In what follows, the multiplication sign is omitted for the sake of brevity. The neutral elements with respect to addition and multiplication coincide with the arithmetic zero $0$ and one $1$.

The notion and notation of inverse elements with respect to multiplication, and of exponents have the usual sense. However, the tropical addition $\oplus$ is not invertible (opposite numbers do not exist), and hence a subtraction operation is undefined.

Matrices and vectors over max-algebra are introduced in the usual way. Matrix operations follow the standard rules, where the arithmetic addition $+$ is replaced by the tropical addition $\oplus$. For matrices $\bm{A}=(a_{ij})$, $\bm{B}=(b_{ij})$, $\bm{C}=(c_{ij})$ of compatible sizes, and scalar $x$, matrix addition, multiplication and multiplication by scalar are given by
\begin{equation*}
(\bm{A}\oplus\bm{B})_{ij}
=
a_{ij}\oplus b_{ij},
\qquad
(\bm{A}\bm{C})_{ij}
=
\bigoplus_{k}a_{ik}c_{kj},
\qquad
(x\bm{A})_{ij}
=
xa_{ij}.
\end{equation*}

Note that multiplication of a matrix (vector) by a scalar gives the same result as in the standard arithmetic. The zero vector denoted by $\bm{0}$, positive vector, zero matrix, and identity matrix denoted by $\bm{I}$ have the standard forms.

In what follows, all vectors are considered column vectors until otherwise indicated. 

For any nonzero column vector $\bm{x}=(x_{j})$, its multiplicative conjugate transpose is the row vector $\bm{x}^{-}=(x_{j}^{-})$, where $x_{j}^{-}=x_{j}^{-1}$ if $x_{j}\ne0$, and $x_{j}^{-}=0$ otherwise. 

Multiplicative conjugate transposition of a nonzero matrix $\bm{A}=(a_{ij})$ results in the transposed matrix $\bm{A}^{-}=(a_{ij}^{-})$, where $a_{ij}^{-}=a_{ji}^{-1}$ if $a_{ji}\ne0$, and $a_{ij}^{-}=0$ otherwise.

For a square matrix $\bm{A}$, the nonnegative integer powers indicate repeated tropical multiplication of $\bm{A}$ by itself, defined as $\bm{A}^{0}=\bm{I}$, $\bm{A}^{p}=\bm{A}^{p-1}\bm{A}$ for all integers $p>0$.

The tropical trace of an $(n\times n)$-matrix $\bm{A}=(a_{ij})$ is given by
\begin{equation*}
\mathop\mathrm{tr}\bm{A}=a_{11}\oplus\cdots\oplus a_{nn}
=
\bigoplus_{i=1}^{n}a_{ii}.
\end{equation*}

The spectral radius of the matrix $\bm{A}$ is calculated by the formula
\begin{equation*}
\lambda
=
\mathop\mathrm{tr}\bm{A}\oplus\cdots\oplus\mathop\mathrm{tr}\nolimits^{1/n}(\bm{A}^{n})
=
\bigoplus_{k=1}^{n}{\mathop\mathrm{tr}}^{1/k}(\bm{A}^{k}).
\end{equation*}

Under the condition that $\lambda\leq1$, the Kleene operator (the Kleene star) is defined to send the matrix $\bm{A}$ into the matrix
\begin{equation*}
\bm{A}^{\ast}
=
\bm{I}\oplus\bm{A}\oplus\cdots\oplus\bm{A}^{n-1}
=
\bigoplus_{k=0}^{n-1}\bm{A}^{k}.
\end{equation*}

For any $(n\times m)$-matrix $\bm{A}=(a_{ij})$ and $n$-vector $\bm{x}=(x_{j})$, the tropical analogues of the matrix and vector norms are given by
\begin{equation*}
\|\bm{A}\|
=
\bigoplus_{i,j=1}^{n}a_{ij},
\qquad
\|\bm{x}\|
=
\bigoplus_{j=1}^{n}x_{j}.
\end{equation*}

\subsection{Tropical Representation of Log-Chebyshev Approximation} 

In this section, we show how the optimization problems involved in the implementation of log-Chebyshev approximation can be represented and solved in the framework of max-algebra where addition is defined as maximum, and multiplication as usual.

We start with problem \eqref{P-minx_maxijlogaijlogxixj}, which can be written in terms of max-algebra in vector form as
\begin{equation}
\begin{aligned}
\min_{\bm{x}>\bm{0}}
&&&
\bm{x}^{-}\bm{A}\bm{x}.
\end{aligned}
\label{P-minxxAx}
\end{equation}

The next result offers a complete solution of problem \eqref{P-minxxAx} (see, e.g., \cite{Krivulin2014Constrained,Krivulin2015Multidimensional}).
\begin{lemma}
\label{L-minxxAx}
Let $\bm{A}$ be a $(n\times n)$-matrix with spectral radius $\lambda>0$. Then, the minimum in problem \eqref{P-minxxAx} is equal to $\lambda$, and all positive solutions are given in the parametric form
\begin{equation*}
\bm{x}
=
(\lambda^{-1}\bm{A})^{\ast}
\bm{u},
\qquad
\bm{u}>\bm{0}.
\end{equation*}
\end{lemma}

Note that the parametric form of the result defines the set of solutions as a tropical linear span of the columns in the corresponding matrix that generates all solutions.

Furthermore, we define an $(n\times n)$-matrix $\bm{B}$ with columns $\bm{b}_{1},\ldots,\bm{b}_{n}$ as follows:
\begin{equation*}
\bm{B}
=
(\lambda^{-1}\bm{A})^{\ast}.
\end{equation*}

We observe that the solution set \eqref{E-X} can be written as
\begin{equation*}
X
=
\left\{
(\lambda^{-1}\bm{A})^{\ast}\bm{u}|\ \bm{u}>\bm{0}
\right\}
=
\{\bm{B}\bm{u}|\ \bm{u}>\bm{0}\}.
\end{equation*}

Consider the maximization problem at \eqref{P-maxxmaxi1xi-maxixi1-minxmaxi1xi-maxixi1}, which now takes the form
\begin{equation}
\begin{aligned}
\max_{\bm{x}\in X}
&&&
\|\bm{x}^{-}\|;
\\
\text{s.t.}
&&&
\|\bm{x}\|
=
1;
\end{aligned}
\label{P-maxxx-xeq1}
\end{equation}

A complete solution of the problem is given as follows \cite{Krivulin2024Application}.
\begin{lemma}
The maximum in problem \eqref{P-maxxx-xeq1} is $\|\bm{B}\bm{B}^{-}\|$, and all solutions are given by
\begin{equation*}
\bm{x}
=
\bm{b}_{k}\|\bm{b}_{k}\|^{-1},
\qquad
k
=
\arg\max_{1\leq j\leq n}\|\bm{b}_{j}\|\|(\bm{b}_{j})^{-}\|.
\end{equation*}
\end{lemma} 

In the context of the pairwise comparison problems, we are interested in the minimal solution of the problem. Therefore, if a vector of the set of solutions is dominated by the other vectors or is the only vector in the set, it uniquely determines the minimal solution of the problem. Otherwise, the minimal solution is not uniquely defined. 

The minimization problem at \eqref{P-maxxmaxi1xi-maxixi1-minxmaxi1xi-maxixi1} is represented as
\begin{equation}
\begin{aligned}
\min_{\bm{x}\in X}
&&&
\|\bm{x}^{-}\|;
\\
\text{s.t.}
&&&
\|\bm{x}\|
=
1;
\end{aligned}
\label{P-minxx-xeq1}
\end{equation}

The solution of the problem is unique and described by the next result \cite{Krivulin2024Application}.
\begin{lemma}
The minimum in problem \eqref{P-minxx-xeq1} is equal to $\|\bm{B}\|$, and the solution is given by
\begin{equation*}
\bm{x}
=
(\bm{1}^{T}\bm{B})^{-}.
\end{equation*}
\end{lemma}

\subsection{Solution of Multicriteria Pairwise Comparison Problem}
\label{S-SMPCP}

We now describe the minimax weighted LCA solution procedure \cite{Krivulin2019Tropical,Krivulin2019Methods,Krivulin2024Application}, which uses results of tropical optimization including the solutions presented in the previous section. The procedure consists of the following steps described in terms of max-algebra.
\begin{enumerate}
\renewcommand{\labelenumii}{\arabic{enumi}.\arabic{enumii}.}
\item
Determination of the best and worst differentiating vectors of weights for the criteria.
\begin{enumerate}
\item
Construction of the generating matrix $\bm{D}$ for the vectors of weights:
\begin{equation}
\bm{D}
=
(\lambda^{-1}\bm{C})^{\ast}
=
\bigoplus_{i=0}^{m-1}
(\lambda^{-1}\bm{C})^{i},
\qquad
\lambda
=
\bigoplus_{i=1}^{m}
\operatorname{tr}^{1/i}(\bm{C}^{i}).
\label{E-D-lambda}
\end{equation}
\item
Calculation of the best differentiating vector of weights $\bm{w}=(w_{i})$ from the matrix $\bm{D}$ with columns $\bm{d}_{1},\ldots,\bm{d}_{m}$ as the componentwise smallest of the vectors 
\begin{equation}
\bm{w}
=
\bm{d}_{k}\|\bm{d}_{k}\|^{-1},
\qquad
k
=
\arg\max_{1\leq i\leq m}\|\bm{d}_{i}\|\|\bm{d}_{i}^{-}\|.
\label{E-w-k}
\end{equation}
If the smallest vector cannot be uniquely determined, then all vectors that are not greater than any other vector are considered best.
\item
Calculation of the worst differentiating vector of weights $\bm{v}=(v_{i})$ given by
\begin{equation}
\bm{v}
=
(\bm{1}^{T}\bm{D})^{-}.
\label{E-v}
\end{equation}
\end{enumerate}
\item
Determination of the best differentiating vector of ratings for alternatives.
\begin{enumerate}
\item
Calculation of the weighted sum of pairwise comparison matrices of alternatives:
\begin{equation}
\bm{P}
=
\bigoplus_{i=1}^{m}
w_{i}\bm{A}_{i}.
\label{E-P}
\end{equation}
\item
Construction of the generating matrix $\bm{Q}$ for the vectors of ratings for alternatives:
\begin{equation}
\bm{Q}
=
(\mu^{-1}\bm{P})^{\ast}
=
\bigoplus_{j=0}^{n-1}
(\mu^{-1}\bm{P})^{j},
\qquad
\mu
=
\bigoplus_{j=1}^{n}
\operatorname{tr}^{1/j}(\bm{P}^{j}).
\label{E-Q-mu}
\end{equation}
\item
Calculation of the best differentiating vector of ratings from the matrix $\bm{Q}$ with columns $\bm{q}_{1},\ldots,\bm{q}_{n}$ as the componentwise smallest of the vectors 
\begin{equation}
\bm{x}
=
\bm{q}_{l}\|\bm{q}_{l}\|^{-1},
\qquad
l
=
\arg\max_{1\leq j\leq n}\|\bm{q}_{j}\|\|\bm{q}_{j}^{-}\|.
\label{E-x-m}
\end{equation}
If the smallest vector cannot be uniquely determined, then all vectors that are not greater than any other vector are considered best.
\end{enumerate}
\item
Determination of the worst differentiating vector of ratings for alternatives.
\begin{enumerate}
\item
Calculation of the weighted sum of pairwise comparison matrices of alternatives:
\begin{equation}
\bm{R}
=
\bigoplus_{i=1}^{m}
v_{i}\bm{A}_{i}.
\label{E-R}
\end{equation}
\item
Construction of the generating matrix $\bm{S}$ for the vectors of ratings for alternatives:
\begin{equation}
\bm{S}
=
(\nu^{-1}\bm{R})^{\ast}
=
\bigoplus_{j=0}^{n-1}
(\nu^{-1}\bm{R})^{j},
\qquad
\nu
=
\bigoplus_{j=1}^{n}
\operatorname{tr}^{1/n}(\bm{R}^{j}).
\label{E-S-nu}
\end{equation}
\item
Calculation the worst differentiating vector of ratings:
\begin{equation}
\bm{y}
=
(\bm{1}^{T}\bm{S})^{-}.
\label{E-y}
\end{equation}
\end{enumerate}
\end{enumerate}

\section{School Selection Example}
\label{S-SSE}

To illustrate the computational technique used, we consider a known problem from \cite{Saaty1977Scaling,Saaty1990Analytic} of rating three high schools A, B and C to select the most preferable one for admission. The schools are compared according to the six criteria: (i)~learning, (ii)~friends, (iii)~school life, (iv)~vocational training, (v)~college preparation, (vi)~music classes. The results of pairwise comparison of the criteria are given by the matrix 
\begin{equation*}
C
=
\begin{pmatrix}
1 & 4 & 3 & 1 & 3 & 4
\\
1/4 & 1 & 7 & 3 & 1/5 & 1
\\
1/3 & 1/7 & 1 & 1/5 & 1/5 & 1/6
\\
1 & 1/3 & 5 & 1 & 1 & 1/3
\\
1/3 & 5 & 5 & 1 & 1 & 3
\\
1/4 & 1 & 6 & 3 & 1/3 & 1
\end{pmatrix}.
\end{equation*}

Pairwise comparisons of the alternatives A, B and C under each criterion yield the following matrices 
\begin{gather*}
A_{1}
=
\begin{pmatrix}
1 & 1/3 & 1/2
\\
3 & 1 & 3
\\
2 & 1/3 & 1
\end{pmatrix},
\quad
A_{2}
=
\begin{pmatrix}
1 & 1 & 1
\\
1 & 1 & 1
\\
1 & 1 & 1
\end{pmatrix},
\quad
A_{3}
=
\begin{pmatrix}
1 & 5 & 1
\\
1/5 & 1 & 1/5
\\
1 & 5 & 1
\end{pmatrix},
\\
A_{4}
=
\begin{pmatrix}
1 & 9 & 7
\\
1/9 & 1 & 1/5
\\
1/7 & 5 & 1
\end{pmatrix},
\quad
A_{5}
=
\begin{pmatrix}
1 & 1/2 & 1
\\
2 & 1 & 2
\\
1 & 1/2 & 1
\end{pmatrix},
\quad
A_{6}
=
\begin{pmatrix}
1 & 6 & 4
\\
1/6 & 1 & 1/3
\\
1/4 & 3 & 1
\end{pmatrix}.
\end{gather*}

Below we present and compare the results of solving the problem using AHP, WGM and LCA methods.

\subsection{Solution by Analytic Hierarchy Process}

The solution using AHP requires finding the principal eigenvectors for all pairwise comparison matrices normalized by the sum of their entries. For the pairwise comparison matrix $\bm{C}$ for criteria, we have the eigenvector
\begin{equation*}
\bm{w}
\approx
\begin{pmatrix}
0.3208 & 0.1395 & 0.0348 & 0.1285 & 0.2374 & 0.1391
\end{pmatrix}^{T}.
\end{equation*}

The eigenvectors of the pairwise comparison matrices $\bm{A}_{1},\ldots,\bm{A}_{6}$ are written as
\begin{gather*}
\bm{x}_{1}
\approx
\begin{pmatrix}
0.1571 \\
0.5936 \\
0.2493
\end{pmatrix},
\qquad
\bm{x}_{2}
\approx
\begin{pmatrix}
0.3333 \\
0.3333 \\
0.3333
\end{pmatrix},
\qquad
\bm{x}_{3}
\approx
\begin{pmatrix}
0.4545 \\
0.0909 \\
0.4545
\end{pmatrix},
\\
\bm{x}_{4}
\approx
\begin{pmatrix}
0.7720 \\
0.0545 \\
0.1734
\end{pmatrix},
\qquad
\bm{x}_{5}
\approx
\begin{pmatrix}
0.2500 \\
0.5000 \\
0.2500
\end{pmatrix},
\qquad
\bm{x}_{6}
\approx
\begin{pmatrix}
0.6910 \\
0.0914 \\
0.2176
\end{pmatrix}.
\end{gather*}

After calculating the vector of ratings for alternatives using formula \eqref{E-AHP}, normalizing it relative to the maximum entry and determining the order of alternatives, we obtain
\begin{equation*}
\bm{x}
\approx
\begin{pmatrix}
0.3673 \\
0.3785 \\
0.2542
\end{pmatrix},
\qquad
\bm{x}/\max_{j}x_{j}
\approx
\begin{pmatrix}
0.9705 \\
1.0000 \\
0.6715
\end{pmatrix},
\qquad
\text{B}\succ\text{A}\succ\text{C}.
\end{equation*}

\subsection{Solution by Weighted Geometric Means}

To apply the method of weighted geometric means, we find the vectors of geometric means for the row entries of the pairwise comparison matrices. The pairwise comparison matrix $\bm{C}$ of criteria has a vector of geometric means normalized with respect to the sum of its entries in the form
\begin{equation*}
\bm{w}
\approx
\begin{pmatrix}
0.3160 & 0.1391 & 0.0360 & 0.1251 & 0.2360 & 0.1477
\end{pmatrix}^{T}.
\end{equation*}

Calculating the vectors of geometric means for the matrices $\bm{A}_{1},\ldots,\bm{A}_{6}$ yields
\begin{gather*}
\bm{x}_{1}
\approx
\begin{pmatrix}
0.5503 \\
2.0801 \\
0.8736
\end{pmatrix},
\qquad
\bm{x}_{2}
\approx
\begin{pmatrix}
1.0000 \\
1.0000 \\
1.0000
\end{pmatrix},
\qquad
\bm{x}_{3}
\approx
\begin{pmatrix}
1.7100 \\
0.3420 \\
1.7100
\end{pmatrix},
\\
\bm{x}_{4}
\approx
\begin{pmatrix}
3.9791 \\
0.2811 \\
0.8939
\end{pmatrix},
\qquad
\bm{x}_{5}
\approx
\begin{pmatrix}
0.7937 \\
1.5874 \\
0.7937
\end{pmatrix},
\qquad
\bm{x}_{6}
\approx
\begin{pmatrix}
2.8845 \\
0.3816 \\
0.9086
\end{pmatrix}.
\end{gather*}

The vector of ratings for alternatives found using formula \eqref{E-WGM}, the result of its normalization by the maximum element and the order of alternatives are written as
\begin{equation*}
\bm{x}
\approx
\begin{pmatrix}
1.1111 \\
1.0008 \\
0.8993
\end{pmatrix},
\qquad
\bm{x}/\max_{j}x_{j}
\approx
\begin{pmatrix}
1.0000 \\
0.9007 \\
0.8094
\end{pmatrix},
\qquad
\text{A}\succ\text{B}\succ\text{C}.
\end{equation*}

\subsection{Solution by Log-Chebyshev Approximation}

This section describes the solution obtained in accordance with the procedure in Section~\ref{S-SMPCP}, where all calculations are performed in terms of the max-algebra, for which addition is defined as max, and multiplications as usual.

We first take the pairwise comparison matrix $\bm{C}$ and apply formulas at \eqref{E-D-lambda} to construct a generating matrix $\bm{D}$ for the vectors of weights of the criteria. As a result of successive calculations, we obtain 
\begin{equation*}
\lambda
\approx
2.5900,
\qquad
\bm{D}
=
(\lambda^{-1}\bm{C})^{\ast}
\approx
\begin{pmatrix}
1.0000 & 2.2361 & 6.0434 & 2.5900 & 1.1583 & 1.5444 \\
0.4472 & 1.0000 & 2.7027 & 1.1583 & 0.5180 & 0.6907 \\
0.1287 & 0.2878 & 1.0000 & 0.3333 & 0.1491 & 0.1988 \\
0.3861 & 0.8633 & 2.3333 & 1.0000 & 0.4472 & 0.5963 \\
0.8633 & 1.9305 & 5.2175 & 2.2361 & 1.0000 & 1.3333 \\
0.4472 & 1.0000 & 2.7027 & 1.1583 & 0.5180 & 1.0000
\end{pmatrix}.
\end{equation*}

Next we use \eqref{E-w-k} and \eqref{E-v} to find the best and worst differentiating vectors of weights 
\begin{equation*}
\bm{w}
\approx
\begin{pmatrix}
1.0000 \\
0.4472 \\
0.1287 \\
0.3861 \\
0.8633 \\
0.4472
\end{pmatrix},
\qquad
\bm{v}
\approx
\begin{pmatrix}
1.0000 \\
0.4472 \\
0.1655 \\
0.3861 \\
0.8633 \\
0.6475
\end{pmatrix}.
\end{equation*}

Determination of the best differentiating vector of ratings for alternatives starts with the calculation of the matrix $\bm{P}$ according to \eqref{E-P}, which yields 
\begin{equation*}
\bm{P}
=
w_{1}\bm{A}_{1}
\oplus
\cdots
\oplus
w_{6}\bm{A}_{6}
\approx
\begin{pmatrix}
1.0000 & 3.4749 & 2.7027 \\
3.0000 & 1.0000 & 3.0000 \\
2.0000 & 1.9305 & 1.0000
\end{pmatrix}.
\end{equation*}

We follow \eqref{E-Q-mu} to form a generating matrix for vectors of ratings. As a result, we have  
\begin{equation*}
\mu
\approx
3.2287,
\qquad
\bm{Q}
=
(\mu^{-1}\bm{P})^{\ast}
\approx
\begin{pmatrix}
1.0000 & 1.0762 & 1.0000 \\
0.9292 & 1.0000 & 0.9292 \\
0.6194 & 0.6667 & 1.0000
\end{pmatrix}.
\end{equation*}

Application of \eqref{E-x-m} leads to a single best differentiating vector of ratings for the alternatives and the corresponding order of alternatives in the form 
\begin{equation*}
\bm{x}
\approx
\begin{pmatrix}
1.0000 \\
0.9292 \\
0.6194
\end{pmatrix},
\qquad
\text{A}\succ\text{B}\succ\text{C}.
\end{equation*}

To determine the worst differentiating vector of ratings for alternatives, we first use \eqref{E-R} to calculate the matrix $\bm{R}$ and write 
\begin{equation*}
\bm{R}
=
v_{1}\bm{A}_{1}
\oplus
\cdots
\oplus
v_{6}\bm{A}_{6}
\approx
\begin{pmatrix}
1.0000 & 3.8850 & 2.7027 \\
3.0000 & 1.0000 & 3.0000 \\
2.0000 & 1.9425 & 1.0000
\end{pmatrix}.
\end{equation*}

The construction of the generating matrix according to \eqref{E-S-nu} results in
\begin{equation*}
\nu
\approx
3.4140
\qquad
\bm{S}
=
(\nu^{-1}\bm{R})^{\ast}
\approx
\begin{pmatrix}
1.0000 & 1.1380 & 1.0000 \\
0.8787 & 1.0000 & 0.8787 \\
0.5858 & 0.6667 & 1.0000
\end{pmatrix}.
\end{equation*}

After calculation of the worst differentiating vector of ratings using \eqref{E-y}, we arrive at the following vector of ratings and the order of alternatives: 
\begin{equation*}
\bm{y}
\approx
\begin{pmatrix}
1.0000 \\
0.8787 \\
1.0000
\end{pmatrix},
\qquad
\text{A}\equiv\text{C}\succ\text{B}.
\end{equation*}

The result obtained demonstrate that the ranks of alternatives based on ratings obtained by different methods may differ. Specifically, the AHP solution assigns the highest rating to alternative $\mathcal{A}_{2}$ (the school B), while the other methods choose alternative $\mathcal{A}_{1}$ (the school A). At the same time, the ranks of the alternatives provided by the WGM solution and the best differentiating LGA solution coincide. 

In the case when the results of different methods conflict with each other, a decision on selecting the best alternative can be made on the basis of collection and comparison of the results of all methods. The alternative that receives the highest rating more often than other alternatives as a result of applying these methods is a very reasonable candidate to be chosen as the best. In accordance with these considerations, the most preferable alternative for the above school selection example should be $\mathcal{A}_{1}$ (the school A) indicated as the best by two of the three solutions considered.

\bibliographystyle{abbrvurl}

\bibliography{Application_of_log-Chebyshev_approximation_and_tropical_algebra_to_multicriteria_problems_of_pairwise_comparisons}

\end{document}